# A Note On The Kadison-Singer Problem

Charles A. Akemann, Betül Tanbay and Ali Ülger

**Abstract**. Let $H$ be a separable Hilbert space with a fixed orthonormal basis $(e_n)_{n\geq 1}$ and $B(H)$ be the full von Neumann algebra of the bounded linear operators $T: H \to H$. Identifying $\ell^\infty = C(\beta N)$ with the diagonal operators, we consider $C(\beta N)$ as a subalgebra of $B(H)$. For each $t \in \beta N$, let $[\delta_t]$ be the *set* of the states of $B(H)$ that extend the Dirac measure $\delta_t$. Our main result shows that, for each $t$ in $\beta N$, the set $[\delta_t]$ either lies in a finite dimensional subspace of $B(H)^*$ or else it must contain a homeomorphic copy of $\beta N$.

**Introduction**. Let $H$ be a separable infinite dimensional Hilbert space with a fixed orthonormal basis $(e_n)_{n\geq 1}$. Let $B(H)$ be the full von Neumann algebra of the bounded linear operators $T: H \to H$. Identifying each bounded sequence with a diagonal operator, we can consider $\ell^\infty$ as a $C^*$–subalgebra of $B(H)$. The famous Kadison-Singer problem asks the following:

*Problem* (**KS**): Does every pure state of $\ell^\infty$ extend in a unique way to a pure state of $B(H)$?

This 50 year old problem has turned out to be a basic problem related to a dozen other important problems [Cz-Tr]. An apparently more general problem, which is extensively studied in the papers [Ar1], [Ar2], [Ar3], [Ch-Ku] and [Bu-Ch], is the following one:

*Problem* (**A**). Let $A \subseteq B$ be two $C^*$-algebras, $A$ being a $C^*$-subalgebra of $B$. When does every pure state of $A$ extend in a unique way to a pure state of $B$?

Problem (KS), and also partially problem (A), may be generalized as follows.

*Problem* (**B**). Let $K$ be a compact Hausdorff space and $Y$ be a Banach space such that $C(K) \subseteq Y$ (i.e. $Y$ contains $C(K)$ as a closed subspace). When does every Dirac measure $\delta_t$ ($t \in K$) extend in a unique way to a functional on $Y$?

Hoping that this general problem may shed some light on the Kadison-Singer problem, we want to study Problem (B) under various hypotheses. Suppose for a moment that, for each $t \in K$, the Dirac measure $\delta_t$ has a unique Hahn-Banach extension to an element $\overline{\delta_t}$ of $Y^*$. Then $\|\overline{\delta_t}\| = 1$ and, as one can easily see (see Lemma 1 in [Ch-Ku]), the mapping $t \mapsto \overline{\delta_t}$ is a weak$^*$ continuous function from $K$ into $Y^*$. Whether $\delta_t$ extends to $Y$ uniquely or not, the set $[\delta_t]$ of all the Hahn-Banach

---

This work of the authors has been supported in part by IMBM, Istanbul Center for Mathematics.

extensions of $\delta_t$ of norm 1 is a weak* compact convex subset of $Y^*$ and, for $t \neq s$, the sets $[\delta_t]$ and $[\delta_s]$ are disjoint. So we can always choose (by the axiom of choice) an element from each of these sets. In this way we define a "selection mapping" $\rho : K \to Y^*$ such that for each $t \in K$, $\rho(t)$ is a norm preserving extension of the functional $\delta_t$ to $Y$. If such a function $\rho$ exists and is weak* continuous, we say that the pair $(K, Y)$ has the "**continuous extension property**". Of course there is no reason why such a continuous $\rho$ would exist. However, for instance, if the multi-valued mapping $t \mapsto [\delta_t]$ is lower semi continuous in the weak* topology of $Y^*$ then, by Michael's Selection Theorem [Mi], such a continuous $\rho$ exists. Also, if the space $Y^*$ has the Kadec-Klee property ( i.e. on the unit sphere $\|y^*\|= 1$ of $Y^*$ the weak-star and the norm topologies agree) then again such a continuous $\rho$ exists.

The first main result of the paper says that the pair $(K, Y)$ has the continuous extension property iff the space $C(K)$ is complemented in the space $Y$ by a contractive projection. Concerning Problem (B), the main result is the following: Suppose that the space $Y^*$ has the property (V) of Pelczynski and $C(K)$ does not contain an isomorphic copy of $\ell^\infty$. Then the pair $(K, Y^*)$ has the continuous extension property iff $K$ is finite. The third main result of the paper is the result stated in the abstract.

**1**-**Preliminary Results and Terminology**. In this section we recall the definitions of the Banach space properties used in the subsequent sections. Throughout this section $X$ and $Y$ will be two arbitrary Banach spaces and $T : X \to Y$ a bounded linear operator. We always regard $X$ as naturally embedded into its second dual $X^{**}$.

*Weakly Unconditionally Cauchy Series.* A series $\sum_{n=0}^{\infty} x_n$ in the Banach space $X$ is said to be wuC if, for each $f \in X^*$, $\sum_{n=0}^{\infty}|f(x_n)|< \infty$.

*Unconditionally Converging Operators.* The operator $T : X \to Y$ is said to be unconditionally converging if it transforms wuC series in $X$ into unconditionally converging series in $Y$.

*Pelczynski's Property (V).* We recall that the space $X$ has property (V) iff any unconditionally converging linear operator from $X$ into any other Banach space is weakly compact [Pe]. Any nonreflexive Banach space having the property (V) contains an isomorphic copy of $c_0$ [Pe].

*Grothendieck Property.* The Banach space $X$ is said to have the Grothendieck property if any weak-star convergent sequence in $X^*$ converges weakly. Grothendieck proved that the space $\ell^\infty$ has this property [Gr]. As proved by Pfitzner [Pf], actually any von Neumann algebra has the property (V), so the Grothendieck property as well. We recall that any dual space having the property (V) has the Grothendieck property [Di].

**2**-**Main Results**. Let $Y$ be an arbitrary Banach space and $K$ be any compact Hausdorff space. Suppose that $Y$ contains an isometric copy of the space $C(K)$ so that we can and do consider $C(K)$ as a subspace of $Y$. Let us recall that the pair $(K, Y)$ is said to have the **continuous extension property** if there is a continuous mapping $\rho$ from $K$ into $(Y^*, w^*)$ such that, for each $t \in K$, $\rho(t)$ is a norm preserving extension of the functional $\delta_t$ to an element of $Y^*$.

The next result gives us some information about the question when the pair $(K, Y)$ has the continuous extension property.

**Lemma 2.1**. Let $Y$ be a Banach space and $K$ be any compact Hausdorff space. Suppose that $Y$ contains $C(K)$ as a closed subspace. Then the pair $(K, Y)$ has the continuous extension property iff there is a contractive projection from $Y$ onto $C(K)$.

*Proof.* Suppose first that the pair $(K, Y)$ has the continuous extension property. So we have a mapping $\rho : K \to Y^*$, which is continuous for the weak-star topology of $Y^*$ and such that, for $a \in C(K)$, $< a, \rho(t) > = < a, \delta_t > = a(t)$. Let $\varphi : Y \to C(K)$ be the mapping defined by

$$\varphi(y)(t) = < y, \rho(t) >.$$

The operator $\varphi$ sends the space $Y$ into $C(K)$ since $\rho$ is continuous for the weak-star topology of $Y$. Moreover $\varphi$ is linear and continuous. The restriction of $\varphi$ to the subspace $C(K)$ of $Y$ is just the identity mapping on $C(K)$. Hence $\varphi$ is a bounded projection from $Y$ onto $C(K)$. As $\|\varphi\| \leq 1$ (actually $\|\varphi\| = 1$ since $\varphi(1_K) = I_K$), the projection $\varphi$ is a contractive projection.

Conversely, let $P : Y \to C(K)$ be a contractive projection. Then its adjoint $P^* : M(K) \to Y^*$ is continuous in the weak-star topologies of the corresponding spaces and, for each $t \in K$, $\|P^*(\delta_t)\| = 1$. Since on the Gelfand spectrum $\{\delta_t : t \in K\}$ of $C(K)$, the weak-star topology induced by $\sigma(M(K), C(K))$ is the same as the original topology of $K$, the mapping $\rho : K \to Y^*$, defined by $\rho(t) = P^*(\delta_t)$, is continuous from $K$ into $(Y^*, w^*)$. Moreover, for $a \in C(K)$,

$$< \rho(t), a > = < P^*(\delta_t), a > = < \delta_t, P(a) > = < \delta_t, a > = a(t)$$

so that $\rho$ is a continuous extension mapping. Hence the pair $(K, Y)$ has the continuous extension property. ∎

Let us recall that a bounded projection $P : Y \to Y$ is said to be an $M$–projection if, for each $y \in Y$, $\|y\| = \max\{\|P(y)\|, \|y - P(y)\|\}$. If $P : Y \to Y$ is an $M$–projection and $Q$ is any contractive projection on $Y$ with $P(Y) = Q(Y)$ then $P = Q$. [HWW; p.2, Proposition 1.2]. From this fact the next result follows immediately.

**Corollary 2.2**. If $C(K)$ is the range of an $M$-projection $P : Y \to Y$ then the mapping $\rho : K \to Y^*$, $\rho(t) = P^*(\delta_t)$, is the only continuous extension map from $K$ into $Y^*$. ∎

To proceed we need the following result.

**Lemma 2.3**. Let $X$ and $Y$ be two Banach spaces. If the space $Y$ does not contain an isomorphic copy of $\ell^\infty$ then every bounded linear operator $T : X^* \to Y$ is unconditionally converging. In particular every bounded linear operator $T : X^* \to Y$ is weakly compact if $X^*$ has the property (V) and $Y$ does not contain an isomorphic copy of $\ell^\infty$.

*Proof.* Suppose that $Y$ does not contain an isomorphic copy of $\ell^\infty$. Let $T : X^* \to Y$ be a bounded linear operator. If $T$ is not unconditionally converging then $X^*$ has a subspace $M$ isomorphic to $c_0$ such that the restriction of $T$ to $M$ is an isomorphism from $M$ onto $T(M)$ [Pe]. Let $i : M \to X$ be the natural injection. Then, since $M^{**}$ is isomorphic to $\ell^\infty$ and $Y$ does not contain an isomorphic copy of $\ell^\infty$, by a result of Rosenthal [Ro; Proposition 1.2], the linear operator
$T^{**} \circ P \circ i : M^{**} \to X^{***} \to X^* \to Y$ is weakly compact. Here $P : X^{***} \to X^*$ is the natural (i.e. restriction) projection. It follows that the restriction of the mapping $T^{**} \circ P \circ i$ to $M$, which is just the restriction of $T$ to $M$, is weakly compact. As the restriction of $T$ to $M$ is an isomorphism, it cannot be weakly compact. This contradiction proves that $T$ is unconditionally converging. The last assertion follows from the characterizations of the spaces having the property (V) given by Pelczynski in [Pe]∎

As proved by Pfitzner [Pf], every von Neumann algebra $B$ has the property (V), hence the Grothendieck property. That is, the weak-star convergent sequences in $B^*$ converge weakly. Actually it is possible to extract from Pfitzner's work [Pf] a considerably stronger result. Apparently, this result has not been previously observed.

**Theorem 2.4**. Let $B$ be a von Neumann algebra and $K$ a weak-star compact subset of $B^*$. Then $K$ is weakly compact iff it does not contain a homeomorphic copy of $\beta N$.

*Proof.* If $K$ is weakly compact then it cannot contain a homeomorphic copy of $\beta N$ since the weakly compact subsets of any Banach space are weakly sequentially compact whereas the space $\beta N$ does not contain any convergent infinite sequence. Conversely, suppose that $K$ does not contain a homeomorphic copy of $\beta N$. Then the space $C(K)$ does not contain an isomorphic copy of $\ell^\infty$. [Ha; p. 67, second Remark]. Now let $\varphi : B \to C(K)$ be the linear operator defined by $\varphi(m)(f) = <m, f>$. Since the space $C(K)$ does not contain an isomorphic copy of $\ell^\infty$, by the preceding lemma, $\varphi$ is unconditionally converging. Hence, since $B$ has the property (V), $\varphi$ is weakly compact. For $f \in K$, let $\delta_f$ be the Dirac measure at $f$. Then $\varphi^*(\delta_f) = f$, so that $K \subseteq \varphi^*(X)$. Here $X$ is the closed unit ball of the Banach space $M(K)$ of the regular Borel measures on $K$, the dual space of $C(K)$. It follows that $K$ is weakly compact.∎

Since $Card(\beta N) = 2^c$ [Wi; p.140] (here $c = Card(R)$), the preceding theorem implies that any weak-star compact subset $K$ of the dual of a von Neumann algebra with $Card(K) < 2^c$ is weakly compact. Thus every net $(f_\alpha)_{\alpha \in I}$ in the dual of a von Neumann algebra $B$ that lies in a weak-star compact subset $K$ of $B^*$ with

$Card(K) < 2^c$ has a weakly convergent subnet. It is clear from the proof of the preceding theorem that, the conclusion of Theorem 2.4 is not special to von Neumann algebras. It is also valid for any dual Banach space that has the property (V).

It is well-known that $c_0$ is not complemented in $\ell^\infty$. As a general version of this result we give the following corollary. This result extends some known results (see e.g. [Ch-Ku; Theorem 6] and [Bu-Ch; Lemma 3.5]).

**Corollary 2.5**. Let $B$ be a von Neumann algebra and $A$ be a $C^*$-subalgebra of $B$ which does not contain an isomorphic copy of $\ell^\infty$. Then there is a bounded projection from $B$ onto $A$ iff the dimension of $A$ is finite.

*Proof.* Let $P : B \to A$ be a bounded projection of $B$ onto $A$. As the algebra $A$ does not contain an isomorphic copy of $\ell^\infty$, by Lemma 2.3 above, $P$ is weakly compact. Then, since $P$ is onto, by the Open Mapping Theorem, the closed unit ball of the algebra $A$ is weakly compact. Hence $A$ is reflexive. In particular $A$ is weakly sequentially complete. Hence, by [Sa; Proposition 2], the dimension of $A$ is finite. The converse is trivial.■

For pairs of $C^*$-algebras $(A, B)$, where $A$ is a $C^*$-subalgebra of $B$, such that there is a unique projection of norm one from $B$ onto $A$, we refer the reader to the papers of Archbold mentioned above.

As a result related to Theorem 2.4 we mention the following result of Anderson. In [An1; Theorem 6], Anderson proves under the Continuum Hypothesis that $\beta N$ has an infinite compact subset $K$ such that, for each $t \in K$, the Dirac measure $\delta_t$ has a unique state (so pure state) extension $\overline{\delta_t}$ to $B(H)$. So the mapping $\rho : K \to B(H)^*$, $\rho(t) = \overline{\delta_t}$, is a continuous extension mapping. By the Tietze Extension Theorem, the restriction mapping from $C(\beta N)$ onto $C(K)$ is a bounded surjective linear operator. As this mapping is onto, so not weakly compact, by Lemma 2.3, the space $C(K)$ contains an isomorphic copy of $\ell^\infty$. Combined with Theorem 2.4, this result of Anderson shows that although the space $C(K)$ contains an isomorphic copy of the space $\ell^\infty$, every pure state of $C(K)$ extends in a unique way to a pure state of $B(H)$. In the opposite direction, again under Continuum Hypothesis, Ch. Akemann and N. Weaver have proved in [Ak-We] that there exists a pure state $f \in B(H)^*$ whose restriction to any masa is not pure.

We also recall the following result. We include a proof for the sake of completeness.

**Lemma 2.6**. Let $X$ and $Y$ be two Banach spaces such that $X$ is a subspace of $Y^*$. If $X$ is complemented in $Y^*$ then $X$ is complemented in its second dual $X^{**}$.

*Proof.* Let $p$ be the natural projection from $Y^{***} = Y^* \oplus Y^\perp$ onto $Y^*$. Let $q : Y^* \to X$ be a bounded projection. We consider $X^{**}$ as naturally embedded into the space $Y^{***}$. Then the composition $p \circ q^{**}$ is a projection that sends $X^{**}$ onto $X$.■

The second main result of this paper is the following result, which is very closely related to the above mentioned result of Anderson.

**Theorem 2.7**. Let $Y$ be an Banach space and $K$ an infinite compact Hausdorff space. Suppose that $Y^*$ contains $C(K)$ as a closed subspace. If the pair $(K, Y^*)$ has the continuous extension property, then the space $C(K)$ is complemented in its second dual. In particular, in this case the space $C(K)$ has the Grothendieck property and contains an isomorphic copy of $\ell^\infty$.

*Proof.* The assertion that $C(K)$ is complemented in its second dual follows directly from the preceding lemma and Lemma 2.1. As $C(K)^{**}$ has the Grothendieck property, any complemented subspace of it, in particular $C(K)$, has the Grothendieck property. If the space $C(K)$ did not contain an isomorphic copy of $\ell^\infty$, any bounded projection from $C(K)^{**}$ onto $C(K)$ would be weakly compact by Lemma 2.3. This is not possible unless $C(K)$ is reflexive, which is not the case since $K$ is infinite. This contradiction proves that $C(K)$ contains an isomorphic copy of $\ell^\infty$.∎

**Remarks 2.8**. **a**) In [Ha] R. Haydon has constructed a compact Hausdorff space $K$ such that the space $C(K)$ has the Grothendieck property and yet the space $C(K)$ does not contain an isomorphic copy of $\ell^\infty$. However if the space $C(K)$ contains an isomorphic copy of $\ell^\infty$ then the compact $K$ contains a homeomorphic copy of $\beta N$ [Ha; p.67, second Remark].

**b**) Even if a compact Hausdorff space $K$ contains a homeomorphic copy of $\beta N$, the space $C(K)$ may not contain an isomorphic copy of $\ell^\infty$ (see [Ha, p. 67] and related reference there). So in the unique extension problem of the pure states of the algebra $C(K)$ to a larger von Neumann algebra, the essential hypothesis is not the existence of a homeomorphic copy of $\beta N$ in $K$ but the existence of an isomorphic copy of $\ell^\infty$ in $C(K)$.

**3**-**Kadison-Singer Problem**. In this section we present a result directly related to the Kadison-Singer problem. To this end, let $H$ be a separable Hilbert space with a fixed orthonormal basis $(e_i)_{i \in N}$. For a bounded sequence $\lambda = (\lambda_n)_{n \in N}$, let $T_\lambda : H \to H$ be the bounded linear operator that sends an element $x = \sum_{n=0}^{\infty} x_i e_i$ of $H$ to the element $T_\lambda(x) = \sum_{n=0}^{\infty} \lambda_i x_i e_i$. The correspondence $\lambda \mapsto T_\lambda$ is an $*$–isometry from $\ell^\infty$ into $B(H)$. Identifying $\ell^\infty$ with its image under this isometry, we can and do consider $\ell^\infty$ as a von Neumann subalgebra of $B(H)$. The mapping $D : B(H) \to \ell^\infty$, defined by $D(T) = (< T(e_n), e_n >)_{n \in N}$, is a contractive positive projection: $P(T^*) = \overline{P(T)}$ and $P(T^*T)$ is a positive sequence in $\ell^\infty$. Further, we identify the space $\ell^\infty$ with the abelian $C^*$–algebra $C(\beta N)$. For each $t$ in $\beta N$, let $\delta_t$ be the Dirac measure at $t$. Then $D^*(\delta_t)$ is given by

$$D^*(\delta_t)(T) = \lim_t \ < T(e_n), e_n >.$$

The right hand side of the preceding equality denotes the limit of the bounded sequence $(< T(e_n), e_n >)_{n \in N}$ over the ultrafilter $t$. Since the projection $D$ is positive and contractive, for each $t \in \beta N$, $D^*(\delta_t)$ is a state (actually pure state [An1]) extension of the Dirac measure $\delta_t$ to $B(H)$ and the function $\rho : \beta N \to B(H)^*$, defined by $\rho(t) = D^*(\delta_t)$, is a *continuous extension mapping*. The projection $D$ induces the decomposition

$$B(H) = \ell^\infty \oplus B_0(H),$$

where $B_0(H)$ is the kernel of $D$. It follows that $B(H)^*$ decomposes as

$$B(H)^* = \ell^{\infty *} \oplus \ell^{\infty \perp},$$

where $\ell^{\infty \perp}$ is the annihilator of $\ell^\infty$ in $B(H)^*$.

Fix now an element $t$ in $\beta N$. Since $B(H)^* = \ell^{\infty *} \oplus \ell^{\infty \perp}$, every extension of $\delta_t$ to $B(H)$ is of the form $\rho = \delta_t + \lambda$. Here $\lambda \in \ell^{\infty \perp}$ so that $\lambda$ vanishes on the diagonal operators. For each $t \in \beta N$, by $[\delta_t]$ we denote the set of all the state extensions of $\delta_t$ to the algebra $B(H)$. The subset $[\delta_t]$ of $B(H)^*$ is *convex and weak-star compact*.

To proceed we need some preliminary results.

Let $A$ be an arbitrary $C^*$–algebra. By $P(A)$ we denote the set of the pure states of $A$. For $\tau \in P(A)$, let $N_\tau = \{a \in A : \tau(a^*a) = 0\}$. The set $N_\tau$ is a maximal modular left ideal and $N_\tau + N_\tau^* = Ker(\tau)$ [Ka], where $N_\tau^* = \{a^* : a \in N_\tau\}$. The ideal $N_\tau$ has a bounded right approximate identity consisting of an increasing net of positive elements in the closed unit ball of $N_\tau$. It follows that the left ideal $N_\tau^{**}$ ( the second dual of $N_\tau$ considered as an ideal in the von Neumann algebra $A^{**}$) has a right unit, denoted $e_\tau$. This $e_\tau$ is a positive idempotent, so a projection, in the von Neumann algebra $A^{**}$. Also, $Ker(\tau)^{**} = \{m \in A^{**} :< m, \tau >= 0\}$. Here $Ker(\tau)^{**}$ denotes the second dual of the Banach space $Ker(\tau)$, which is identified with the weak-star closure of $Ker(\tau)$ in $A^{**}$. Let $1$ denotes the unit element of the von Neumann algebra $A^{**}$. Since $e_\tau \in N_\tau^{**}$, we have $< e_\tau, \tau >= 0$. Moreover, since $e_\tau$ is self-adjoint and a right unit in $N_\tau^{**}$, for $m \in N_\tau^{**}$, we have $m.e_\tau = m$ and $e_t.m^* = m^*$. These facts will be used in the proof of the next lemmas. See also [Ped, 3.3.16] and [Ch-Ku; Lemma 5].

**Lemma 3. 1**. For each $\tau \in P(A)$, there is a minimal projection $e \in A^{**}$ such that $< \tau, e >= 1$. If for some other pure state $\tau^{'}$ we have $< \tau^{'}, e >= 1$ then $\tau = \tau^{'}$.

*Proof.* We fix a pure state $\tau \in P(A)$. Then, with the above notation, $< \tau, 1 - e_\tau >= 1$. Our first aim is to prove that the projection $e = 1 - e_\tau$ is a minimal projection in $A^{**}$. To this end, first observe that, for $a \in N_\tau$, we have $eae = ea(1 - e_\tau) = 0$. Similarly, for $a \in N_\tau$, we have $ea^*e = (1 - e_\tau)a^*e = 0$. Hence, for each $a \in Ker(\tau)$, we have $eae = 0$. Since the multiplication in $A^{**}$ is separately weak-star continuous, by the weak-star density of $Ker(\tau)$ into $Ker(\tau)^{**}$, we get that for all $m \in Ker(\tau)^{**}$, we have $eme = 0$. Since $A^{**} = Ker(\tau)^{**} \oplus Ce$, we conclude that $e$ is a minimal projection in $A^{**}$.

Now for $a \in A$,

$$|<\tau, ae_\tau>|^2 \leq <\tau, e_\tau>.\tau(a^*a) = 0$$

so that $<\tau, ae_\tau> = 0$. Hence $<\tau, ae> = <\tau, a(1-e_\tau)> = <\tau, a>$. Similarly, $<\tau, ea> = <\tau, a>$. Thus, for all $a \in A$, $<\tau, eae> = <\tau, ea> = <\tau, a>$. Since $<\tau', e> = 1$ by hypothesis, we also have $<\tau', eae> = <\tau', a>$. As $e$ is a minimal projection in $A^{**}$, for each $a \in A$, we have $eae = \lambda e$ for some constant $\lambda$. Hence, for all $a \in A$,

$$<\tau, a> = <\tau, eae> = \lambda <\tau, e> = \lambda = \lambda <\tau', e> = <\tau', eae> = <\tau', a>$$

so that $\tau = \tau'$. ∎

**Lemma 3.2**. Let $B$ be a unital $C^*$-algebra and $B$ a $C^*$-subalgebra of $B$ sharing the same unit. Consider $A^{**}$ as a von Neumann subalgebra of $B^{**}$, and let $\tau$ be a pure state of $A$. Then, in the notation of the paragraph preceding lemma 3.1, the set of all states of $B$ that restrict to $\tau$ on $A$ is exactly the set of states $f$ of $B$ such that $f(1 - e_\tau) = 1$. Further, this set spans a finite dimensional subspace of $B^*$ iff the space $(1 - e_\tau)B^{**}(1 - e_\tau)$ is finite dimensional. Finally, the pure state $\tau$ has unique state extension to $B$ iff the space $(1 - e_\tau)B^{**}(1 - e_\tau)$ is one dimensional (so that $(1 - e_\tau)B^*(1 - e_\tau)$ contains a unique state).

*Proof.* Since $(1 - e_\tau)A^{**}(1 - e_\tau)$ is one dimensional (as shown above), then clearly $f(1 - e_\tau) = 1$ implies that $f|_A = \tau$. The converse is also immediate, and the rest of the lemma is also clear. ∎

The third main result of the paper is a dichotomy theorem that classify the sets $[\delta_t]$ as "very large" and "small".

**Theorem 3.3**. For each $t \in \beta N$, either the set $[\delta_t]$ lies in a finite dimensional subspace of $B(H)^*$ or it contains a homeomorphic copy of $\beta N$.

*Proof.* Using Theorem 2.4, it will suffice to show that the set $[\delta_t]$ is not weakly compact if it does not lie in a finite dimensional subspace of $B(H)^*$. So suppose that the set $[\delta_t]$ is not weakly compact. Using the notation developed in the previous two lemmas, we need only show that the set $[\delta_t]$ contains a sequence of pure states $\tau_n$ that are orthogonal in the sense that their supporting projections $e'_{\tau_n}$ s in $B(H)^{**}$ are orthogonal. Since the set $[\delta_t]$ does not lie in a finite dimensional subspace of $B(H)^*$, by the last lemma, the space $(1 - e_\tau)B(H)^*(1 - e_\tau)$ is infinite dimensional. Thus the set of states in that subspace spans an infinite dimensional subspace also. Since that set is weak* closed (see Proposition 3.11.9 in [Ped]), the Krein-Milman theorem implies that the set of pure states in $(1 - e_\tau)B(H)^*(1 - e_\tau)$ must also span an infinite dimensional space. Since each pure state is supported by a minimal projection, this means that the set of minimal projections spans an infinite dimensional subspace of $(1 - e_\tau)B(H)^{**}(1 - e_\tau)$. Consequently there is a sequence of orthogonal minimal projections in $(1 - e_\tau)B(H)^{**}(1 - e_\tau)$, and each

such minimal projection supports a pure state $\tau_n$. The theorem follows.∎

We here remark that in the case where $[\delta_t]$ lies in a finite dimensional subspace of $B(H)^*$ it is norm compact, so norm separable. After this theorem the following question becomes crucial.

**Question**. How to prove that, for a given $t \in \beta N$, the set $[\delta_t]$ is weakly compact; or equivalently, does not contain a homeomorphic copy of $\beta N$?

Now let $t \in \beta N$ be a given ultrafilter and $\rho = \overline{\delta_t} + \lambda$ be a pure state of $B(H)$ extending $\delta_t$. Then the set

$$N_\rho = \{T \in B(H) : \rho(T^*T) = 0\}$$

is the closed maximal left ideal associated to the pure state $\rho$. The ideal $N_\rho$ is in general neither weak-star closed nor has a right unit. But it always has a positive bounded right approximate identity $(U_i)_{i \in I}$. For certain $\delta_t$, $(t \in \beta N)$, as this is the case for $\delta_n$ $(n \in N)$, $N_\rho$ may have a bounded right approximate identity consisting of positive diagonal operators. So it is not unreasonable to expect that, for certain $\rho's$, each $U_i$ is a diagonal operator. Actually we have the following result, which follows directly from lemma 3.2 above.

**Theorem 3.4**. Let $t \in \beta N$ be a given point, and let $\rho = \overline{\delta_t} + \lambda$ be a pure state extension of $\delta_t$ to $B(H)$. Then the maximal left ideal $N_\rho$ has a positive bounded right approximate identity $(U_i)_{i \in I}$ consisting of diagonal operators iff the set $[\delta_t]$ is a singleton so that $\rho = \overline{\delta_t}$ is the unique pure state extension of $\delta_t$ to $B(H)$ in this case.∎

Next we want to study some topological properties of the union of the sets $[\delta_t]$.

Let $E$ be the set of those $t \in \beta N$ such that the set $[\delta_t]$ is weakly (so norm) compact. The set $E$ is nonempty since it contains the set of integers. If the Continuum Hypothesis is assumed, then the set $E$ contains much more than the integers. Let $\sum = \cup_{t \in E} [\delta_t]$ be the union of all the sets $[\delta_t]$ for $t \in E$. On the set $\sum$ we put the metric induced by the norm of $B(H)^*$. We denote this metric by $d$.

**Lemma 3.5**. The metric space $(\sum, d)$ is complete and locally compact. Moreover, for each subset $F$ of $E$, the set $\sum' = \cup_{t \in F} [\delta_t]$ is both open and closed in $\sum$.

*Proof.* For $t \neq s$ ($t, s$ in $\beta N$), there exists (by Urysohn Lemma) a diagonal operator $T$ such that $\|T\| = 1$, $< \delta_t, T > = 1$ and $< \delta_s, T > = -1$. So, for any $\rho \in [\delta_t]$ and $\rho' \in [\delta_s]$, $\|\rho - \rho'\| \geq |< \delta_t + \lambda - \delta_s - \lambda', T >| = 2$. Hence, the metric distance $d([\delta_t], [\delta_s])$ between the closed sets $[\delta_t]$ and $[\delta_s]$ is 2. It follows that the set $\sum$ is norm closed in the space $B(H)^*$. So the metric space $(\sum, d)$ is complete. Moreover, for the same reasons, for any nonempty subset $F$ of $E$, the set $\sum' = \cup_{t \in F} [\delta_t]$ is

also closed in $B(H)^*$, so in $\Sigma$. Since the complement of $\Sigma'$ in $\Sigma$ is also open in $\Sigma$, the set $\Sigma'$ is both open and closed in $\Sigma$. In particular, for each $t \in E$, the set $[\delta_t]$ is both compact and open in $\Sigma$. This in turn shows that each $\rho \in \Sigma$ has a compact neighborhood so that the metric space $(\Sigma, d)$ is locally compact.∎

Let $\Sigma$ be as in the preceding lemma. Since the metric space $(\Sigma, d)$ is locally compact, we can consider its Stone-Cech compactification $\beta\Sigma$. Take a $t \in \beta N$ for which the set $[\delta_t]$ does not contain a homeomorphic copy of $\beta N$. Then the set $[\delta_t]$ is a compact-open subset of $\Sigma$. So it is also compact and open in the space $\beta\Sigma$. Hence the characteristic function of each set $[\delta_t]$ is in the space $C(\beta\Sigma)$. Thus the space $C([\delta_t])$ is a complemented ideal of the $C^*$-algebra $C(\beta\Sigma)$, so that we have a bounded projection $P : C(\beta\Sigma) \to C([\delta_t])$. The compact $[\delta_t]$ being metric, the $C^*$-algebra $C([\delta_t])$ is separable. So if the space $C(\beta\Sigma)$ has the Grothendieck property the projection $P$ is weakly compact. This implies that the $C^*$-algebra $C([\delta_t])$ is finite dimensional, which in turn implies that the set $[\delta_t]$ is finite. As this set is convex, this is possible only if it contains a single point. Thus, if the space $C(\beta\Sigma)$ has the Grothendieck property, then the sets $[\delta_t]$ are one-point sets whenever they do not contain $\beta N$ homeomorphically. Whence the questions

**Questions**. 1-Does the space $C(\beta\Sigma)$ have the Grothendieck property ?
               2-Is the set $E$ closed in $\beta N$?

14 (2) (1964), 659-664.

[Wi] S. Willard, General Topology, Addison Wesley Publishing Company (1970) Reading, Massachusets.

Department of Mathematics, University of California, Santa Barbara, Ca 93106.

Department of Mathematics, Bogazici University , Bebek-Istanbul, Turkey.

Department of Mathematics, Koc University , 34450 Sariyer-Istanbul,Turkey.

e-mails:

    akemann@math.ucsb.edu

    tanbay@boun.edu.tr

    aulger@ku.edu.tr